\documentclass[11pt,a4paper]{amsart}
\usepackage{amsfonts,amssymb,bbm}
\usepackage{enumerate,mathrsfs}\usepackage{graphicx}

\def\k{\mathbbm{k}}

\def\N{\Bbb{N}}

\def\R{\Bbb{R}}

\def\di{\partial}

\def\liml{\lim\limits}
\def\suml{\sum\limits}\def\intl{\int\limits}

\def\capl{\mathop\cap\limits}

\newcommand{\quots}[2]{{\footnotesize\left.\raisebox{0.4ex}{$#1$}\! / \!\raisebox{-0.4ex}{$#2$}\right.}}

\def\tb{\tilde{b}}

\def\tF{\tilde{F}}

\def\tN{{\tilde{N}}}

\def\ty{\tilde{y}}
\def\tuy{\tilde{\uy}}

\def\hf{\hat{f}}

\def\hg{{\hat{g}}}

\def\hR{{\widehat{R}}}

\def\hy{\hat{y}}
\def\huy{{\hat{\underline{y}}}}

\def\al{\alpha}\def\de{\delta}\def\De{\Delta}
\def\ep{\epsilon}

\def\ca{\mathfrak a}
\def\cb{\mathfrak b}

\def\cm{{\frak m}}

\def\cq{{\frak q}}
\def\cU{\mathcal U}
\def\cV{\mathcal V}

\def\uc{{\underline{c}}}

\def\uk{{\underline{k}}}\def\ul{{\underline{l}}}\def\um{{\underline{m}}}
\def\ux{{\underline{x}}}\def\uy{{\underline{y}}}\def\uz{{\underline{z}}}

\def\one{{1\hspace{-0.1cm}\rm I}}

\newcommand{\ber}{\begin{array}{l}}\newcommand{\eer}{\end{array}}
\newcommand{\bpm}{\begin{pmatrix}}\newcommand{\epm}{\end{pmatrix}}
\newcommand{\bbm}{\begin{bmatrix}}\newcommand{\ebm}{\end{bmatrix}}
\newcommand{\bM}{\begin{matrix}}\newcommand{\eM}{\end{matrix}}
\newcommand{\bee}{\begin{enumerate}}\newcommand{\eee}{\end{enumerate}}
\newcommand{\bei}{\begin{itemize}}\newcommand{\eei}{\end{itemize}}

\def\sset{\subset}\def\sseteq{\subseteq}\def\ssetneq{\subsetneq}\def\smin{\setminus}

\newtheorem{Lemma}{Lemma}[section]\newcommand{\bel}{\begin{Lemma}}\newcommand{\eel}{\end{Lemma}}
\newtheorem{Theorem}[Lemma]{Theorem}\newcommand{\bthe}{\begin{Theorem}}\newcommand{\ethe}{\end{Theorem}}
\newtheorem{Proposition}[Lemma]{Proposition}\newcommand{\bprop}{\begin{Proposition}}\newcommand{\eprop}{\end{Proposition}}
\newtheorem{Corollary}[Lemma]{Corollary}\newcommand{\bcor}{\begin{Corollary}}\newcommand{\ecor}{\end{Corollary}}

\newtheorem{Definition}[Lemma]{Definition}
\newcommand{\bed}{\begin{Definition}}
\newcommand{\eed}{\end{Definition}}
\newtheorem{Definition-Proposition}[Lemma]{Definition-Proposition}

\def\bpr{~\\{\em Proof.\ \ }}
\newcommand{\epr}{{\hfill\ensuremath\blacksquare}}
\newtheorem{Remark}[Lemma]{Remark}
\newcommand{\beR}{\begin{Remark}\rm}
\newcommand{\eeR}{\end{Remark}}
\newtheorem{Example}[Lemma]{Example}
\newcommand{\bex}{\begin{Example}\rm}
\newcommand{\eex}{\end{Example}}
\newtheorem{Problem}[Lemma]{Problem}
\newcommand{\bprob}{\begin{Problem}\rm}
\newcommand{\eprob}{\end{Problem}}

\newcommand{\bet}{\begin{tabular}{cccccccc}}\newcommand{\eet}{\end{tabular}}
\newcommand{\beq}{\begin{equation}}\newcommand{\eeq}{\end{equation}}

\newcommand{\bin}[2]{\binom{#1}{#2}}

\title[]{A\MakeLowercase{pproximation results of} A\MakeLowercase{rtin-}T\MakeLowercase{ougeron-type for general filtrations and for
 }$C^r$-\MakeLowercase{equations.}}
\author[]{G\MakeLowercase{enrich} B\MakeLowercase{elitskii,}
A\MakeLowercase{lberto} F. B\MakeLowercase{oix and}
D\MakeLowercase{mitry} K\MakeLowercase{erner}}
\address{Department of Mathematics, Ben Gurion University of the Negev, P.O.B. 653, Be'er Sheva 84105, Israel.}
\email{genrich@math.bgu.ac.il} \email{fernanal@post.bgu.ac.il}
\email{dmitry.kerner@gmail.com}
\date{\today
}
\thanks{A.F. Boix was partially supported by Israel Science Foundation (grant No. 1910/18) and Spanish Ministerio de Econom\'ia y Competitividad
 MTM2016-7881-P}
\thanks{D.Kerner was partially supported by Israel Science Foundation (grant No. 1910/18)}

\begin{document}
\begin{abstract}
Artin approximation and other related approximation results are used
in various areas. The traditional
 formulation of such results is restricted to filtrations by powers of ideals, $\{I^j\}$, and to Noetherian rings.
 In this paper we extend several approximation results both to rather general filtrations and to $C^r$-rings,
  for  $2\le r\le\infty$.

  As an auxiliary step we establish the surjectivity of the completion map $\quots{C^\infty(\cU)}{J}\to \widehat{\quots{C^\infty(\cU)}{J}}$
   for a very broad class of filtrations.
\end{abstract}

\maketitle
\setcounter{secnumdepth}{6} \setcounter{tocdepth}{1}

All the rings are commutative, unital.
 For a descending filtration by ideals, $R=I_0\supseteq I_1\supseteq\cdots$,  we denote $I_\infty:=\cap I_j$.
We use the multivariable notations, $\ux=(x_1,\dots,x_m)$,
$\uy=(y_1,\dots,y_n)$.

\section{Introduction}
Various versions of Artin approximation are widely used in Algebraic/Analytic Geometry, Commutative Algebra and Singularity Theory.
 Recently they became important in other areas, see  \cite{Rond} for the general introduction and the review of the current state of research.

Traditionally, the approximation statements were restricted to Noetherian rings and to filtrations   by  powers of  ideals, $\{I^j\}$.
 (Two notable exceptions being \cite{Schoutens} and \cite{Moret-Bailly.12}.)

For various recent applications in Singularity Theory
  one needs these approximations both for rings of differentiable/smooth functions
   and for more general filtrations/completions, see
\cite{B.K.motor} and \cite{Boi.Gre.Ker}.
In this note we  extend some of the classical approximation results both to rather general filtrations and to
$C^r$-rings, where $2\leq r\leq\infty.$
This allows, e.g. immediate applications of Artin approximation to the study of non-isolated singularities of maps and schemes.

Below we recall some classical results.

\subsection{Polynomial equations}
Consider a (finite) system of polynomial equations, $F(\uy)\!=\!0$,  where $F(\uy)\!\!\in\!\! R[\uy]^{s}$.
\bed\label{Def.AP}
 The Artin approximation property, AP,
holds for  $R,I_\bullet$ if  for every finite system of polynomial
equations over $R$, a (formal) solution in the completion $\hR^{(I_\bullet)}$
implies the existence of a solution in $R$, which can be chosen
arbitrary close to the formal solution in the filtration topology.
\eed

The famous characterization of  rings with AP reads:
\bthe\label{Thm.Artin.Popescu}\cite[Remark 2.15]{Popescu}, see also \cite[Theorem 1.3]{Popescu.1986} and \cite[Theorem 1]{Rotthaus}

 Let $R$ be a commutative Noetherian excellent ring.
\bee[1.]
\item
If the pair $(R,I)$ is Henselian, for some ideal $I\sset R$, then  AP holds for $R$ and the filtration $I^\bullet$.
\item
If a local ring $(R,\cm)$ has AP, for the filtration $\cm^\bullet$, then it is  Henselian.
 \eee
\ethe

\subsection{Analytic/algebraic equations}
 When the equations $F(\ux,\uy)=0$ are non-polynomial, the formal solution does not imply an  ordinary solution.
 Yet the approximation holds for analytic  equations, more generally for equations given by a $W$-system.
 (A $W$-system is a  Weierstrass system of rings over $\k$, see \cite[page2]{Denef-Lipshitz}.)

\bthe\label{Thm.AP.Artin.Schemmel.Analytic.Algebraic}
 Let $\huy(\ux)\in \k[[\ux]]^n$ be a formal solution, i.e. $ F (\ux,\huy(\ux))=0$, assume $\huy(o)=0$.
\bee[1.]
\item \cite[Theorem 1.1]{Denef-Lipshitz}
 Let $\k$ be either a field or a discrete valuation ring, and suppose the system
of equations $F(\ux,\uy)=0$ is given by a $W$-system, i.e., $ F(
\ux, \uy)\in \k\lceil\!\!\lceil \ux, \uy \rceil\!\!\rceil^s.$
 For every $N\in\N$ there exists  a $W$-solution $\uy(\ux)\in  \k\lceil\!\!\lceil\ux\rceil\!\!\rceil^n$
 satisfying: $\uy(\ux)-\huy(\ux)\in \cm^{N+1}\cdot \k[[\ux]]^n$.
\item
\cite[Theorem 1.2]{Artin.68}, \cite[page 135, Theorem 1]{Wavrik.75}, \cite{Schemmel.1982}
 Let $\k$ be a valued field of arbitrary characteristic,
and suppose that the completion of $\k$ with respect to its absolute value is separable over $\k$.
 Suppose the system of equations $F(x,y)=0$ is $\k$-analytic.
 $ F (\ux,\uy)\in \k\{\ux,\uy\}^s$.
 For every $N\!\in\!\N$ there exists an analytic  solution
  $\uy(\ux)\!\in\! \k\{\ux\}^n$ satisfying: $\uy(\ux)\!-\!\huy(\ux)\!\in\! \cm^{N+1}\!\cdot  \k[[\ux]]^n$.
\eee
 \ethe
We recall the widely used particular cases of this theorem:
\bee[1.]
\item (for $W$-systems) Algebraic equations, i.e. $ F (\ux,\uy)\in \k\langle \ux,\uy\rangle^s$, then part one ensures the approximation by
 an algebraic solution,   $\uy(\ux)\in \k\langle \ux,\uy\rangle^n$.
 \item
(for valued fields) The completion of $\k$ with respect to its absolute value is separable over $\k$, e.g. in the following cases:
 when $\k$ is complete, when $\k$ is perfect, and when $\k$ is discrete, see  \cite[pages 38--39]{Abhyankar-van der Put}).
 Then part two ensures the approximation by $\k$-analytic solution, $\uy(\ux)\in \k\{\ux\}^n$.
\eee

\

For $\R$-analytic equations the approximation statement is much stronger:
\bthe\label{Thm.AP.Tougeron.classical}
 \cite[Theoreme 1.2]{Tougeron1976} Let $F(\ux,\uy)\in \R\{\ux,\uy\}^s$ and assume $\huy_0$ is a formal solution.
  Then there exists a solution $\uy(\ux)\in C^\infty(\R^m,o)^n$, whose Taylor series is  $\huy_0(\ux)$.
\bee[1.]
\item
 Moreover, for any $N\in\N$ there exists an analytic solution,  $\uy_{ann}(\ux)\in \R\{\ux\}^n$, that  is $\cm^N$-homotopic to $\uy(\ux)$.
\item If, moreover, $ F (\ux,\uy)\in \R\langle\ux,\uy\rangle^s$ (algebraic power series)
 then for any $N\in\N$ the approximating solution can be chosen algebraic,
 $\uy(\ux)\in \R\langle\ux\rangle^n$. If in addition $\huy(\ux)\in \R\{\ux\}^n$, then the $\cm^N$-homotopy can be chosen analytic.
  \eee
\ethe

Recall that two solutions, $\uy_0(\ux)$, $\uy_1(\ux)$, are called
$\ca$-homotopic, for an ideal $\ca\sset R$, if there exists a
($C^\infty$/analytic)-family of solutions, $\uy(\ux,t)$,  such that:
\beq
\uy_0(\ux)=\uy(\ux,0),\quad
\uy_1(\ux)=\uy(\ux,1), \quad
 \text{ and } \quad \uy(\ux,t)-\uy_0(\ux)\in I\cdot  C^\infty(\R^m,o)^n \quad \text{ for any $t$.}
 \eeq

\subsection{$C^\infty$-equations}
Let $ F \in C^\infty(\R^m\times\R^n,o)^s$, the ring of smooth
function-germs at the origin  $o\in \R^m\times\R^n$.
 A formal solution of the equation $F (\ux,\uy)=0$ is a power series $\huy_0\in \R[[\ux]]^n$ satisfying $ F (\ux,\huy_0)\in (\ux)^\infty\cdot (C^\infty(\R^m,o))^s$.
This condition is understood in the following sense.
  Borel's lemma ensures the surjectivity of the completion map, $C^\infty(\R^m,o)\twoheadrightarrow \R[[\ux]]$. Thus one takes a(ny)
    Borel-representative $\tuy_0\in C^\infty(\R^m,o)^n$ of $\huy_0$ and verifies $ F (\ux,\tuy_0)\in (\ux)^\infty\cdot (C^\infty(\R^m,o))^s$.
This does not depend on the choice of Borel-representative.

\

The naive generalization of  theorems \ref{Thm.Artin.Popescu}, \ref{Thm.AP.Artin.Schemmel.Analytic.Algebraic},
 \ref{Thm.AP.Tougeron.classical} to
$C^\infty$-equations fails, even for linear equations with
$C^\infty$-coefficients.
\bex
\bee[i.]
\item Take a flat function   $\tau\in (x)^\infty \sset C^\infty(\R^1,o)$, e.g.
$\tau(x)=\Big\{\ber e^{-\frac{1}{x^2}},\   x\neq0 \\0,\  x=0.\eer$.
 Consider the equation $ {\tau}^2(x)y = {\tau}(x) $.
 Every formal power series $\hy\in\R[[x]]$ is a formal solution, but the equation has no  continuous solutions.
 \item
 Take a flat function   $\tau\in (x_1,x_2)^\infty \sset C^\infty(\R^2,o)$ and consider the equation $x_1\cdot y=\tau(x)$.
 Assume  $\tau \not\in (x_1)$, e.g. $\tau$ vanishes only at the origin.
Then  $y=0$ is a formal solution, but there are no continuous solutions.

Here in the first case $\di_y F(x,y)|_{y=0}$ is flat. In the second case $\di_y F(x,y)|_{y=0}$ is not flat, but (considered as a morphism of free modules)
 is far from being surjective, its cokernel is supported on the whole line $\{x_1=0\}$.
\item
More generally, suppose a filtration $I_\bullet$ of some ring $R$ satisfies:
 $I_\infty:=\cap I_j\neq (0)$. Then AP does not hold for $R$, $I_\bullet$.
 For example,   the equation $0\cdot y=b$, where $0\neq b\in I_\infty$,   has many formal solutions, in the sense that $0\cdot y-b\in I_\infty$,
  but has no ordinary solution.
\eee
\eex
 Yet, under   additional assumptions, some approximation results are possible in the $C^\infty$-case.
\bthe\label{Thm.Approximation.Cinfty.known}
\bee[1.]
  \item\cite[\S3.2.2]{van der Put}
Given a  set of polynomials in  one variable with smooth coefficients, $F(\uy)\in
\big(C^{\infty}(\R^1,o)[\uy]\big)^s$, let $A\sseteq
C^{\infty}(\R^1,o)$ be the subalgebra generated by
 the coefficients of $F(\uy)$. Suppose $A\cap \cm^\infty=\{0\}\sset C^{\infty}(\R^1,o)$. Then any formal solution
   $\hat{\uy}_0\in \R[[x]]^n$   lifts to an ordinary solution,
     $\uy_0\in \big(C^{\infty}(\R^1,o)\big)^n$, such that  $F(\uy_0)=0\in (C^{\infty}(\R^1,o))^s$ and $\hat{\uy}_0$ is the Taylor expansion of $\uy_0$.
\item \cite[Theorem 5.3]{B.K.IFT}\label{Thm.Approximation.Cinfty.equations}
Let $ F (\ux,\uy)\in \big(C^\infty(\R^m\times\R^n,o)\big)^s$ and
 suppose the equation $ F (\ux,\uy)=0$  has a formal solution,  $\huy_0(x)$.
 Denote $h_{\huy_0}(\ux):=\det\Big[\frac{\di  F }{\di \uy}\big|_{(\ux,\huy_0(\ux))}\cdot \Big(\frac{\di  F }{\di \uy}\big|_{(\ux,\huy_0(\ux))}\!\Big)^T\Big]$ and suppose
  $h_{\huy_0}\cdot\cm^\infty=\cm^\infty$.
 Then $\huy_0(\ux)$ lifts to an ordinary solution, $\uy_0 \in C^\infty(\R^m,o)^n$, $ F (\ux,\uy_0)=0$,
  whose Taylor series at the origin is  $\huy_0(\ux)$.
  \eee
\ethe
In part 2 the $C^\infty$-function $h_{\huy_0}$ is constructed by taking a Borel representative $\tilde\uy_0\in
\big(C^{\infty}(\R^m,o)\big)^n$ of $\huy_0$.
  As before, the condition
 $h_{\huy_0}\cdot\cm^\infty=\cm^\infty$
   does not depend on the choice of the representative.

\subsection{Our results}
\bei
\item In \S\ref{Sec.AP.general.filtrations}
 we reduce the verification  of AP for $R$, $I_\bullet$ to AP for $R$, $I^\bullet$, under very weak assumptions on $I_j$.
  In particular, this extends part 1  of Theorem \ref{Thm.Artin.Popescu} to rather general filtrations $I_\bullet$.
   Similarly we extend theorem \ref{Thm.AP.Artin.Schemmel.Analytic.Algebraic}.

Theorem \ref{Thm.AP.Tougeron.classical} is extended to the general filtration in \S\ref{Sec.AP.Tougeron.Generalization}.

The importance of these results is clear: finer filtrations ensure finer approximations.
\item
In \S\ref{Sec.Approximation.C.infty.equations} we extend part 2. of
theorem \ref{Thm.Approximation.Cinfty.known} to the ring
$\quots{C^\infty(\R^p,o)}{J}$ and the general  filtration $I_\bullet$.
  Moreover, we strengthen it, in the spirit  of theorem \ref{Thm.AP.Tougeron.classical},
  to ensure a solution that is analytic/algebraic modulo the ideal of flat functions, $I_\infty$.

In this section we assume the surjectivity of the completion map  $\quots{C^\infty(\R^p,o)}{J}\to \widehat{\quots{C^\infty(\R^p,o)}{J}}$.
 For general filtrations this question is more complicated than the classical (Borel) surjectivity $C^\infty(\R^p,o) \twoheadrightarrow \R[[\ux]]$.
 We obtain the sufficient condition for the surjectivity   in the appendix.

\item In \S\ref{Sec.Approximation.Cr.equations} we extend part 2 of theorem \ref{Thm.Approximation.Cinfty.known} to $C^r$-equations.
\eei

\subsection{Acknowledgement} We thank M.Sodin for the highly useful reference to \cite{Hormander} and the referee for helpful remarks.

\section{Artin-type approximation for general filtrations}\label{Sec.AP.general.filtrations}
\subsection{The case of polynomial equations}
Let $R$ be  a commutative  (not necessarily Noetherian)  ring, with a filtration $\{I_\bullet\}$. The following condition is a weakening of
 being finitely generated:
\beq\label{Eq.cond.weak.f.g}
\text{for any $N$ there exists $\tN=\tN(N)\gg1$ and a finite set $\{q_\al\}$ in $I_N$ such that $I_{N+\tN}\sseteq (\{q_\al\})$}.
\eeq
\bel\label{Thm.Approximation.Artin.general.filtration}
 Suppose $R$ has AP   for   a filtration $I_\bullet$.
  Then $R$ has AP for any   filtration $ \ca_\bullet$ satisfying condition \eqref{Eq.cond.weak.f.g}
   and such that $\ca_{d_j}\sseteq I_j$, for any $j$ and a corresponding $d_j<\infty$.
\eel
\bpr Let $F(y)\in R[y]^s$ be a system of polynomial equations.
 We should prove: any $\hR^{(\ca_\bullet)}$-formal solution is $\ca_\bullet$-approximated by a solution in $R$.

Take the completion
$R\stackrel{\phi^{(\ca_\bullet)}}{\to}\hR^{(\ca_\bullet)}$ and let
$\hy_0\in (\hR^{(\ca_\bullet)})^n$ be a formal solution.
 For any $N$ and any $\tilde{N}\gg N$ there exists $y_N\in R^n$ (not necessarily a solution) such that $\hy_0-\phi(y_N)\in \ca_{\tilde{N}+1}\cdot (\hR^{(\ca_\bullet)})^n$.

By the assumption \eqref{Eq.cond.weak.f.g} there exists a finite  set of elements $\{q_\al\}\in \ca_{N+1}$   such that
 $\hy_0-\phi(y_N)\in (\{q_\al\})\cdot (\hR^{(\ca_\bullet)})^n$.
Change the variable, $y=y_N+\sum \ty_\al q_\al$. The initial system
of
 equations becomes $F(y_N+\sum \ty_\al q_\al)=0$,   for the unknowns $\{\ty_\al\}$. This system has a $\hR^{(\ca_\bullet)}$-formal solution, coming from $\hy_0$.

By the assumption $\ca_{d_j}\sseteq I_j$,   we have the natural map
$\hR^{(\ca_\bullet)}\stackrel{\psi}{\to} \hR^{(I_\bullet)}$. (It is  not necessarily injective.)
 This map sends the $\hR^{(\ca_\bullet)}$-formal solution to a  $\hR^{(I_\bullet)}$-formal solution:
\beq \phi^{(I_\bullet)}(F) \big( \psi(\hy_0)
\big)=\psi\phi^{(\ca_\bullet)}(F)(\psi(\hy_0))=\psi(0)=0\in
(\hR^{(I_\bullet)})^s.
 \eeq

Now, by AP for $ I_\bullet $-filtration, we get an ordinary solution,
$F(y_N+\sum \ty_\al q_\al)=0$, for some $\{\ty_\al\in R^n\}$. But
then $y_N+\sum \ty_\al q_\al\in R^n$
 is the needed ordinary solution. (It approximates $\hy_0$ for the filtration $\ca_\bullet$.)
\epr
\bex\label{Ex.AP.general.filtrations}
\bee[i.]
\item
Suppose two filtrations are equivalent, $I_\bullet\sim \ca_\bullet$, then $R$ has AP for $I_\bullet$ iff it has AP for $ \ca_\bullet$.
\item  For a Noetherian local ring, $(R,\cm)$, many filtrations satisfy $\cap I_j=0$.
 In particular, for any $j$ and a corresponding $d_j<\infty$ the inclusion  $I_{d_j}\sseteq \cm^j$ holds. Thus AP for $\{\cm^j\}$ implies AP for $I_\bullet$.

\item For the study of non-isolated singularities one needs filtrations of the form $\{\cm^j\cdot J\}$, where the ideal $J$ defines the singular locus.
 (In particular $J$ is not $\cm$-primary.) More generally, one needs filtrations of the form $\{\big(\cap\cq^{n_\al(j)}_\al\big)\cap J\}_j$,
  where $\{\cq_\al\}$ is a finite set of ideals and $\{\liml_{j\to\infty}n_\al(j)= \infty\}_\al$ and $height(J)<height(\cq_\al)$, for any $\al$.
   These filtrations are not equivalent to $ I^\bullet$
   for any $I\sset R$.
   Thus theorem \ref{Thm.Artin.Popescu} cannot be applied directly, but lemma \ref{Thm.Approximation.Artin.general.filtration} is applicable.
\eee
\eex

\subsection{Analytic/$W$-system equations over $\k$}
Theorem \ref{Thm.AP.Artin.Schemmel.Analytic.Algebraic} was   stated for the filtration $ \cm^\bullet$.
 Let  $R$ be one of
   $\quots{\k\{\ux\}}{J}$, $\quots{\k\lceil\!\!\lceil\ux\rceil\!\!\rceil}{J}$.
 (Here $\k$ is a field or a discrete valuation ring, with the assumptions as in theorem \ref{Thm.AP.Artin.Schemmel.Analytic.Algebraic}.)
  Let $F(\ux,\uy)=0$ be the     corresponding system of  equations, i.e. $F\in R\{\uy\}$ or $R\lceil\!\!\lceil\uy\rceil\!\!\rceil$.
\bel\label{Thm.AP.analytic.algebraic.generalization}
Suppose  a  filtration $I_\bullet$ of $R$, satisfies:  $\cm^{d_j}\supseteq I_j$,  for any $j$ and a corresponding $d_j<\infty$.
Suppose the equation $F(\ux,\uy)=0$
 has a formal solution, $\huy_0\in \big(\hR^{(I_\bullet)}\big)^n$.
 For every $N\in\N$ there exists an analytic/$W$-system  solution $\uy_0\in R^n$ satisfying: $\uy_0-\huy_0\in I_{N+1}\cdot (\hR^{(I_\bullet)})^n$.
\eel
 The proof goes by the same argument as in  lemma \ref{Thm.Approximation.Artin.general.filtration}.

\subsection{Analytic  equations over $\R$, a generalization of Tougeron's theorem}\label{Sec.AP.Tougeron.Generalization}
Take the ring $R= \quots{\R\{\ux\}}{J}$  filtered by $I_\bullet$,
 and $F(\ux,\uy)\in (R\{\uy\})^s$. Suppose the equation $F(\ux,\uy)=0$ has a formal solution, $\huy_0\in (\hR^{(I_\bullet)})^n$.

\bprop\label{Thm.AP.Tougeron.Generalization}
\bee[1.]
\item
For any $N\in\N$ there exists  a solution $\uy_0 \in (\quots{C^\infty(\R^m,o)}{J})^n$,
 that  satisfies:
\[
\uy_0-\huy_0\in I_N\cdot\cm^\infty\cdot (\quots{C^\infty(\R^m,o)}{J})^n .
\]
\item
 Moreover, for any $j\in\N$ there exists an analytic solution,  $\uy_{ann}\in R^n$  that  is $I_N\cdot \cm^j$-homotopic to $\uy_0$.
\item If  moreover, $J$ is  generated by algebraic power series and  $ F (\ux,\uy)$ is a (vector of) algebraic power series
 then for any $j\in\N$ the approximating solution can be chosen algebraic,
 $\uy_{alg}\in (\quots{\R\langle\ux\rangle}{J})^n$. If in addition $\huy_0\in (\quots{\R\{\ux\}}{J})^n$,
  then the $I_N\cdot \cm^j$-homotopy can be chosen analytic.
  \eee
\eprop
The condition $\uy_0\!-\!\huy_0\!\in\! I_N\cdot\cm^\infty$ is understood as before: it holds for a(ny) $C^\infty$\!-representative of $\huy_0$.
\bpr
\bee[\bf Step 1.]
\item We reduce the statement to the case $R=\R\{\ux\}$. Let $\tF(\ux,\uy)\in \R\{\ux,\uy\}$ be a representative of $F(\ux,\uy)$.
Fix some (finite) set of generators, $\{q_\al\}$,  of $J$. Consider the equation
\beq\label{Eq.inside.Tougeron}
\tF(\ux,\uy)=\sum_\al q_\al z_\al.
\eeq
Here $\{z_\al\}$ are $s$-columns of new variables. A formal solution of $F(\ux,\uy)=0$ implies a formal solution of \eqref{Eq.inside.Tougeron}.
 Thus, assuming a needed (analytic/algebraic) solution, $\tilde\uy_a$ of \eqref{Eq.inside.Tougeron} (homotopic to the formal solution),
  we get the needed (analytic/algebraic) solution $\uy_a$ of $F(\ux,\uy)=0$, homotopic to $\huy_0$.
\item
Let $R=\R\{\ux\}$ and $F(\ux,\uy)\in \R\{\ux,\uy\}$. Denote by $\uc_0\in R^n$ the $N$'th approximation to the formal solution
 $\huy_0\in (\hR^{(I_\bullet)})^n$, i.e. $\uc_0-\huy_0\in I_N\cdot (\hR)^n$. Fix some generators $\{q_\al\}$ of $I_N$ and consider the shifted equation,
\beq\label{Eq.inside.Tougeron.2}
 F(\ux,\uc_0+\sum_\al q_\al \uz_\al)=0.
\eeq
This is an analytic equation on the new ($n$-columns of) variables $\{\uz_\al\}$. The formal solution $\huy_0$ ensures
 a formal solution $\{\hat\uz_{0,\al}\}$ of
 \eqref{Eq.inside.Tougeron.2}. Then theorem \ref{Thm.AP.Tougeron.classical} ensures $C^\infty$-solutions, $\{ \uz_{0,\al}\}$, whose Taylor
   series are  $\{\hat\uz_{0,\al}\}$.

 Define $\uy_0:=\uc_0+\sum_\al q_\al \uz_{0,\al}\in C^\infty(\R^m,o)^n$. Then $\uy_0-\huy_0\in I_N\cdot \cm^\infty\cdot C^\infty(\R^m,o)^n$.
  Moreover, for any $j\in \N$, Tougeron's theorem ensures analytic solutions, $\{ \uz_{j,\al}\}$ in $\R\{\ux\}$, which are $\cm^j$-homotopic
    to $\{\hat\uz_{0,\al}\}$. This homotopy gives the needed $I_N\cdot \cm^j$-homotopy of $\uy_0$ to $\uy:=\uc_0+\sum q_\al \uz_{j,\al}$.

    This proves parts 1. and 2. of the theorem.

    Part 3. follows similarly, from the $F(\ux,\uy)\in \R\langle\ux\rangle^s$- part of Tougeron's theorem.
\epr\eee
\beR
This proposition is a weak generalization of Tougeron's theorem. One would like to replace the conclusion ``$\uy_0-\huy_0\in I_N\cdot \cm^\infty$" by
 the stronger conclusion $\uy_0-\huy_0\in I_\infty$, i.e. ``$\huy_0$ is the image of $\uy_0$ under the $I_\bullet$-completion". However, this cannot hold without further assumptions.
  Indeed, this would imply (trivially) the surjectivity of the completion map, $C^\infty(\R^m,o)\twoheadrightarrow \widehat{C^\infty(\R^m,o)}\ ^{(I_\bullet)}$.
   But already this surjectivity does not always hold, as it places significant restrictions on the filtration $I_\bullet$..
\eeR

\section{Approximation  for $C^\infty$-equations}\label{Sec.Approximation.C.infty.equations}
 Let $R=\quots{C^\infty(\R^m,o)}{J}$, with some filtration $I_\bullet$.
In this section we always assume  the completion map is  surjective, $R\twoheadrightarrow \hR$.
 This holds for many filtrations, the sufficient condition is established in \ref{Sec.Surjectivity.of.Completion}.
  In particular, the surjectivity holds for filtrations satisfying:
  \beq
(Z,o):=V(I_\infty)=V(I_N),\ \quad \ for\ N\gg1,\quad\quad\quad \{I_{N_j}\sseteq I(Z,o)^j\}_j,\ \text{for\ some }N_j<\infty.
  \eeq
Here $I(Z,o)$ is the ideal of all function-germs that vanish on $(Z,o)$.

\subsection{Formal solutions}

  We often compare elements of $\hR^{(I_\bullet)}$ and $R$.
To simplify the expressions we often put these elements in one formula.
\bee[i.]
\item For $y_1\in R$ and $\hy_0\in \hR^{(I_\bullet)}$ the notation $y_1-\hy_0\in I_j$ means:
 for some representative $y_0\in R$ of $\hy_0$ the difference is $y_1-y_0\in I_j$. (This does not depend on the choice of representative.)

Similarly, the homotopy notation $\hy_0\stackrel{I_j}{\sim}y_1$ means: $y_0\stackrel{I_j}{\sim}y_1$,
 where $y_0$ is a representative of $\hy_0$.
 \item
  For $F(\ux,\uy)\in C^\infty(\R^m\times\R^n,o)$ and $\huy_0\in \hR^n$ the notation $F(\ux,\huy_0)\in I_\infty$ means:
  for
   some representative $\uy_0\in R^n$ of $\huy_0$ the composition satisfies $F(\ux,\uy_0)\in I_\infty$. (This does not depend on the choice of representative.)
 \eee

\

  Take a system of equations, $F (\ux,\uy)=0$, where $F\in \Big(\frac{C^\infty(\R^m\times\R^n,o)}{J}\Big)^s$.

\bed\label{Def.Formal.Solution.Cinfty.equation}
A formal solution is an element   $\huy_0\in  \hR^n$ such
  that $F(\ux,\huy_0)\in  I_\infty\cdot  R^s$.
\eed

\subsection{The approximation theorem}  Suppose there exist a formal solution  $\huy_0\in   \hR^n$.
  Define the auxiliary function-germ as the determinant of the matrix,
\beq\label{Eq.def.of.h}
h_{\huy_0} (\ux):=\det \Big[\frac{\di  F}{\di \uy}\big|_{ (\ux,\huy_0(\ux))}\cdot \Big(\frac{\di  F}{\di \uy}\big|_{ (\ux,\huy_0(\ux))}\Big)^T\Big].
\eeq
As before, in $ F (\ux,\huy_0)$  we substitute   a(ny) $C^\infty$-representative of $\huy_0$. As before, the non-uniqueness of the representative
 changes $h_{\huy_0}  $  only by an element of $I_\infty$.
 The matrix $\frac{\di  F}{\di \uy}\big|_{ (\ux,\huy_0(\ux))}$ is of size $s\times n$, thus $h=0$ unless $n\ge s$.

  We say that $h_{\huy_0}$ has a finite Taylor order at a point $x_0$ if $h_{\huy_0}\not\in\cm^\infty_{x_0}$.

\bthe\label{Thm.Approximation.Cinfty.new}
 Suppose the completion map is surjective, $R\twoheadrightarrow \hR^{(I_\bullet)}$.
Suppose   there exists a formal solution, $ \huy_0 \in  \hR^n $, $\huy_0(0)=0$, satisfying the condition  $h_{\huy_0}\cdot  I_\infty=I_\infty$.
\bee[1.]
\item  There  exists an ordinary solution   $\uy \!\in \! R^n$, such that $ F (\ux,\uy (\ux))\!=\!0$ and
   the $I_\bullet$\!-completion map sends  $\uy $  to $\huy_0$.
 \item Suppose  $F(\ux,\uy)\in \Big(\frac{\R\{\ux,\uy\}+I_\infty\cdot C^\infty(\R^m\times\R^n,o)}{J}\Big)^s$ and the following conditions hold:
 \bee[a.]
 \item the ideals  $J$ and all $I_\bullet$ are analytically generated;
\item   $(Z,o):=V(I_\infty)=V(I_N)$ for $N\gg1$, and $I_\infty\sseteq I(Z,o)^\infty$;
\item  $h_{\huy_0}$ has finite Taylor-orders at all points of $(Z,o)$.
\eee
  Then for any $N\in \N$ exists a solution
  \[
  \uy_N\in \big(\quots{\R\{\ux\}}{J}+I_\infty\big)^n,\quad\quad
     F(\ux,\uy_N(\ux))=0,
\quad\quad    \text{ such that }\uy_N\stackrel{I_N}{\sim}\huy_0.
\]
     \item Suppose
      $F(\ux,\uy)\in \Big(\frac{\R\langle \ux,\uy\rangle+I_\infty\cdot C^\infty(\R^m\times\R^n,o)}{J}\Big)^s$
 and the following conditions hold:
 \bee[a.]
 \item  the ideals $J $ and all $I_\bullet$ are  generated by algebraic power series;
\item       $(Z,o):=V(I_\infty)=V(I_N)$ for $N\gg1$ and  $I_\infty\sseteq I(Z,o)^\infty$;
\item  $h_{\huy_0}$ has finite Taylor-orders at all points of $(Z,o)$.
\eee
  Then for any $N\in \N$ exists a solution
\[
  \uy_N\in \big(\quots{\R\langle\ux\rangle}{J}+I_\infty\big)^n,
  \quad\quad
   F(\ux,\uy_N(\ux))=0
   \quad\quad    \text{ such that }
   \uy_N\stackrel{I_N}{\sim}\huy_0.
\]
 \eee
\ethe
\bpr
\bee[\bf 1.]
\item (The proof expands the initial idea from \cite{B.K.IFT}.)
Let $\tuy\in R^n$ be a $C^\infty$-representative of $\huy_0$, thus $F(x,\tuy)\in I_\infty\cdot R^s$.
 Shift the variables, $\uy  =\tilde\uy  +\De\uy $, and take the Taylor expansion $ F (\ux,\tilde\uy  +\De\uy )$   with remainder:
\beq\label{Eq.expansion.of.f}
 F (\ux,\tuy +\De\uy ) =  F (\ux,\tilde\uy ) + \frac{\di  F (\ux,\tilde\uy )}{\di \uy}\cdot \De\uy +
(\De\uy )^T\left( \int\limits_0^1(1-\xi)\frac
{\partial^2F(\ux,\tilde\uy +\xi\De\uy )}{\partial \uy^2}d\xi\right)(\De\uy ).
\eeq
Thus $ F (\ux,\tilde\uy  +\De\uy ) =0$ is an implicit function equation on $\De\uy$.

We are looking for the solution in the form
\beq
\De\uy(\ux) =h (\ux)\cdot \big( \frac{\di  F (\ux,\tilde\uy )}{\di \uy}\big)^T \cdot
\Big[\frac{\di  F (\ux,\tilde\uy )}{\di \uy}\cdot \Big(\frac{\di  F (\ux,\tilde\uy )}{\di \uy}\Big)^T\Big]^\vee\cdot \uz
\eeq
Here $[\dots]^\vee$ is the adjugate matrix, while $\uz\in R^s$ is a column of free variables.

This substitution gives the equation:
\beq\label{Eq.eventual.IFeq}
\frac{ F (\ux,\tilde\uy )}{h (x)^2}+\uz+\uz^T\cdot \Big[\dots\Big]\uz=0.
\eeq
These are $s$ equations in $s$ variables.

By the assumption $\frac{ F (\ux,\tilde\uy )}{h (x)^2}\in I_\infty\cdot R^s$. The entries of the matrix $ \Big[\dots\Big]$ belong to
 $R$ and depend on $\uz$ via $\De\uy$. Thus they are well defined for any $\uz\in \R^n$, and not just for small values of $\uz$.

Finally, invoke the implicit function theorem in the ring $R$ to get a solution
 $\uz(\ux)\in I_\infty\cdot R^s$.
 This gives the solution  $\uy (\ux)=\tilde\uy (\ux)+\De\uy (\ux)\in R^n$  to $F(\ux,\uy)=0$.

 Note that   $\uy (\ux)$ is sent to $\huy_0(\ux)$ by the completion map,  as was claimed.
 \item
 \bee[\bf Step 1.]
\item   Present $ F = F_{ann}+ F _{flat}$, where $ F_{ann}\in (\quots{\R\{\ux,\uy\}}{J})^s$ and
     $ F _{flat}\in I_\infty\cdot \Big(\quots{C^\infty(\R^m\times\R^n,o)}{J}\Big)^s$.
     If $F(\ux,\huy_0)\in I_\infty\cdot R^s$,       then also $F_{ann}(\ux,\huy_0)\in I_\infty\cdot R^s$. Thus,   the Taylor expansion of $\huy_0$ satisfies:
       $ F_{ann}(\ux,\huy_0)=0$.
      Thus, by proposition  \ref{Thm.AP.Tougeron.Generalization} there exists a family  $\uy(t)\in (\quots{C^\infty((\R^m,o)\times[0,1])}{J})^n$ satisfying:
\beq\ber
\forall\ t:\ \uy(t)-\huy_0\in I_N\cdot\cm^j\cdot R^n,\quad\quad   F_{ann}(\ux,\uy(t))=0,
\\ \uy(0)-\huy_0\in I_N\cdot\cm^\infty\cdot R^n,\quad \quad \uy(1)\in (\quots{\R\{\ux\}}{J})^n.
\eer
\eeq

\item

We verify for any $t$: $h_{\uy(t)}\cdot I_\infty=I_\infty$.
    Indeed, $h_{\huy_0}\cdot I_\infty=I_\infty$ and $h_{\huy_0}$ has finite Taylor-orders at all points of $Z$.
 As $Z$ is closed, and we work with the germ $(Z,o)$, we can assume $Z$ is compact, then this order is bounded.
 Thus there exists a $C^\infty$-representative $\tilde\uy_0$ of $\huy_0$ satisfying for some $d\in \N$:
\beq
h_{\tilde\uy_0}^{-1}(0)\sseteq(Z,o),\quad\quad  \forall\ z\in Z:\ ord_z(h_{\tilde\uy_0})\le d.
\eeq
Thus, for $N\gg1$ and any $t\in [0,1]$ we have: $h_{ \uy(t) }^{-1}(0)\sseteq(Z,o)$, and for any  $z\in Z$:  $ord_z(h_{ \uy(t)})\le d.$
 This implies, for any $t$:  $h_{\uy(t)}\cdot I_\infty=I_\infty$.
\item Finally we consider the equation $F(\ux,\uy(t)+\De(t))=0$, where $\De(t)$ is a (column of) new variable. Expand it as
 in equation \eqref{Eq.expansion.of.f} to get the solution, $\De(t)\in I_\infty\cdot \Big(\quots{C^\infty((\R^m,o)\times[0,1])}{J}\Big)^n$.
  Define $\uy_N:=\uy(1)+\De(1)$, this is a solution, analytic mod $I_\infty$. And our construction ensures $\uy_N\stackrel{I_N}{\sim}\huy_0$.
\eee
\item The proof is the same,  just we use the algebraic part of  proposition \ref{Thm.AP.Tougeron.Generalization}.
\epr \eee

\bex Let $R= \quots{C^\infty(\R^m,o)}{J}$ with a filtration $I_\bullet$ satisfying: $I_\infty\sseteq\cm^\infty$, $V(I_\infty)=V(\cm)=o\in \R^m$.
 This ensures the surjectivity of completion, $R\twoheadrightarrow \hR^{(I_\bullet)}$, see theorem \ref{Thm.Completion.Cinfty.surjective}.
\bee[i.]
\item Suppose  the linear part of the equations is non-degenerate at $0$, i.e. the matrix
  $\frac{\di  F }{\di \uy}\big|_{\substack{\ux=0\\\uy=0}}$ is of rank $s$, with $s\le n$. Then    $h_{\huy_0}$ is invertible for any formal solution $\huy_0$.
   In particular $h_{\huy_0}\cdot\cm^\infty=\cm^\infty$. Thus any formal solution extends to a $C^\infty$-solution.

\item More generally,  assume the derivative $\frac{\di  F }{\di \uy}\big|_{\uy=0}$ is non-degenerate for $\ux\neq0$.
 Thus $h_{\huy_0=0}$ vanishes at $o$ only. Then $h_{\huy_0}\cdot\cm^\infty=\cm^\infty$  holds e.g. if
    $h_{\huy_0=0}$ is analytic.
 This gives a Tougeron type statement for the classical $\cm$-adic completion.
 For $J=(0)$ this gives part 2 of theorem \ref{Thm.Approximation.Cinfty.known}.
\eee
\eex
\bex
Let $R= \quots{C^\infty(\R^m,o)}{J}$ and assume $(Z,o):=V(I_\infty)$ is an analytic germ and moreover: $I_\infty\sseteq I(Z,o)^\infty$,
 and $(Z,o)=V(I_N)$  for $N\gg1$.
By theorem \ref{Thm.Completion.Cinfty.surjective} the completion is surjective again. Given a system of equations, $F(\ux,\uy)=0$, with a formal solution, $\huy_0$, we should check $h_{\huy_0} \cdot I_\infty=I_\infty$.
 Suppose
 $h_{\huy_0}$ is presentable in the form $h_{ann}+h_\infty$, where  $h_\infty\in I_\infty$ and $h_{ann}\in \quots{\R\{\ux\}}{J}$, $h_{ann}^{-1}(0)=Z$.
 (Here we choose some $C^\infty$-representative $\uy_0$ of $\huy_0$, and $h_{ann}$ does not depend on this choice.)

 Then, by {\L}ojasiewicz inequality,  there exist constants $C>0$ and $\de>0$
such that
\[
\text{ $h_{ann}(x)\ge C\cdot dist(x,Z)^\de$ holds in a neighborhood of $(Z,o)$.}
\]
 Therefore  $h_{ann}\cdot I_\infty=I_\infty$ and thus $h_{\huy_0}\cdot I_\infty=I_\infty$.
  Thus theorem \ref{Thm.Approximation.Cinfty.new} ensures a  $C^\infty$-solution, $F(\ux,\uy_0)=0$, whose $I_\bullet$-completion is $\huy_0$.
\eex
\beR
In parts 2,3 of theorem \ref{Thm.Approximation.Cinfty.new} we assume that $F(\ux,\uy)$ is analytic/algebraic modulo $I_\infty$-terms in $\ux$.
 We can allow also the flat terms in $\uy$, i.e.
 ``$F(\ux,\uy)\in \Big(\frac{\R\{\ux,\uy\}+(I_\infty+(\uy)^\infty)\cdot C^\infty(\R^m\times\R^n,o)}{J}\Big)^s$", provided the filtration $I_\bullet$
  satisfies $  \capl_j (I_N)^j\sseteq I_\infty$ for $N\gg1$.
 (The proof goes as before.) This later condition is satisfied for many filtrations.
 \eeR

\section{Approximation  for $C^r$-equations}\label{Sec.Approximation.Cr.equations}
Take the ring of function-germs $C^{r_m,r_n}(\R^m\times\R^n,o)$.
 For $F(\ux,\uy)\in C^{r_m,r_n}(\R^m\times\R^n,o) $
 all the derivatives $\frac{\di^{r_m+r_n} F}{\di x_{i_1}\dots \di x_{i_{r_m}}\di y_{j_1}\dots \di y_{j_{r_n}}}$ exist and are continuous.
 Here $2\le r_m\le r_n \le\infty$. Moreover, if $r_n<\infty$ then we assume $r_m+2\le r_n$.

Fix an ideal $J\sset C^{r_m}(\R^m,o)$ and take the quotient rings,
$\frac{C^{r_m,r_n}(\R^m\times\R^n,o)}{J}$ and $R:=\quots{C^{r_m}(\R^m,o)}{J}$.

An element  $ F \in
\Big(\frac{C^{r_m,r_n}(\R^m\times\R^n,o)}{J}\Big)^s$,
 defines the system of equations, $ F (\ux,\uy)=0$.

\bed
A {\em solution $mod(I)$} to the system $F(\ux,\uy)=0$ is an element  $\tuy_0\in R^n$ satisfying
 $ F (\ux,\tuy_0)\in  I \cdot  R^s$.
\eed
As in the $C^\infty$-case, equation \eqref{Eq.def.of.h}, we define the determinant of the matrix:
\beq
h_{\tuy_0}(\ux):=\det \Big[\frac{\di  F }{\di \uy}\big|_{(\ux,\tuy_0(\ux))}\cdot \Big(\frac{\di  F }{\di \uy}\big|_{(\ux,\tuy_0(\ux))}\Big)^T\Big].
\eeq
 The matrix $\frac{\di  F }{\di \uy}$ is of size $s\times n$, thus $h=0$ unless $n\ge s$.
 The entries of the matrix $\frac{\di  F }{\di \uy}$ lie in $\frac{C^{r_m,r_n-1}(\R^m\times\R^n,o)}{J}$.
 Therefore (as $r_n>r_m$) the entries of the matrix $\frac{\di  F }{\di \uy}\big|_{(\ux,\tuy_0(\ux))}$
 lie in  $R$.

\bprop\label{Thm.Approximation.Cr.new}
 Suppose  $ \tuy_0\in R^n $ is    a $mod(I)$-solution to the equation $ F (\ux,\uy(\ux))=0$, and there holds:
    $I\!\sseteq\! (h_{\tuy_0})^2\!\sset\! R$.
 Then there exists an ordinary solution,  $\uy_0\!\in\! R^n$, such that $ F (\ux,\uy_0(\ux))\!=\!0$ and
   $\uy_0\!-\!\tuy_0\!\in\! \frac{1}{(h_{\tuy_0})^2}I\cdot R^n$.
\eprop
\bpr
 The proof is the same as for theorem \ref{Thm.Approximation.Cinfty.new}.
Shift the variables, $\uy =\tuy_0 +\De\uy$, to get the Taylor expansion as in equation \eqref{Eq.expansion.of.f}

Note that the entries of $\frac{\di  F (\ux,\tilde\uy_0)}{\di \uy}$ and of $\frac
{\partial^2F(\ux,\tilde\uy_0+\xi\De\uy)}{\partial \uy^2}$ belong to $R$, as $r_n\ge r_m+2$.
Thus $F (\ux,\tuy_0 +\De\uy) =0$ is a $C^{r_m}$-implicit function equation.

Proceed as in the proof of theorem \ref{Thm.Approximation.Cinfty.new} to get to equation \eqref{Eq.eventual.IFeq}.

By the assumption $\frac{ F (\ux,\tilde\uy_0)}{h_{ \huy_0}(x)^2}\in \frac{I}{h_{ \huy_0}(x)^2}\cdot  R^s$. The entries of the matrix $ \Big[\dots\Big]$ belong to
 $R$ and depend on $\uz$ via $\De\uy$. Thus they are well defined for any $\uz\in \R^s$, and not just for small values of $\uz$.

Finally, invoke the implicit function theorem in the ring $R$ to get a solution
 $\uz\in I\cdot R^s$.
 This gives the solution  $\uy_0(\ux)=\tilde\uy_0(\ux)+\De\uy(\ux)\in R^n$  to $F(\ux,\uy)=0$.
 Note that   $\uy_0(\ux)$ approximates the initial $\tilde\uy_0(\ux)$, as was claimed.
\epr

\beR
 The assumption $I\sseteq (h_{\tuy_0})^2$ can be weakened.
 Take the annihilator of cokernel of the matrix, $Ann.Coker\big[\frac{\di  F (\ux,\tuy_0)}{\di \uy}\big]\sset R$,
  \cite[\S20]{Eisenbud}. This ideal satisfies
\[
Ann.Coker\big[\frac{\di  F }{\di \uy}\big|_{(\ux,\tuy_0(\ux))}\big]\supseteq
\Big( \det \Big[\frac{\di  F  }{\di \uy}\big|_{(\ux,\tuy_0(\ux))}\cdot \Big(\frac{\di  F  }{\di \uy}\big|_{(\ux,\tuy_0(\ux))}\Big)^T\Big]\Big)
\]
and the proper inclusion often holds.
 Then theorem  \ref{Thm.Approximation.Cr.new} holds with $h$ replaced by any   $\tilde{h}\in Ann.Coker(\frac{\di  F }{\di \uy}\big|_{(\ux,\tuy_0)})$.
 The proof goes with just one change, one puts $\De\uy=\tilde{h}\cdot B\cdot \uz$, where the matrix $B$ satisfies:
  $\frac{\di  F}{\di \uy}\big|_{(\ux,\tuy_0)}\cdot B=\tilde{h}\cdot R^s$.

In the $C^\infty$-case this made no significant difference, as $I$ was the ideal of flat functions. But for $I$ non-flat, the condition $(\tilde{h})^2\supseteq I$
 is often weaker than $(h)^2\supseteq I$.
\eeR

\appendix
\section{The sufficient condition for the surjectivity  of the completion $R\to \hR$}\label{Sec.Surjectivity.of.Completion}
A classical lemma of Borel reads: any real sequence $\{a_{\uk}\}_{\uk\in \N^m}$ is realizable as the sequence of partial derivatives (at $o\in \R^m$) of a
 function $f\in C^\infty(\R^m)$. To specify all the derivatives at $o$ means to specify the Taylor series, thus this lemma
   means the surjectivity of the completion map
   \beq
   C^\infty(\R^m)\twoheadrightarrow \widehat{C^\infty(\R^m)}^{ (\ux)^\bullet}=\R[[\ux]].
   \eeq
 More generally, Whitney's extension theorem
  gives the necessary and sufficient conditions to extend a function with prescribed derivatives on a closed subset $Z\sset \R^m$ to a smooth function on $\R^m$,
    see e.g. \S1.5 of \cite{Narasimhan}.
 In the particular case, $Z$ is a manifold  and the filtration is $\{I_j=I(Z)^j\}$, specifying derivatives on $Z$ is equivalent
  to specifying an element of the completion $\widehat{C^\infty(\R^n)}^{(I_\bullet)}$.
   Then Whitney's extension theorem implies the surjectivity of the completion map
   $C^\infty(\R^n)\to \widehat{C^\infty(\R^n)}^{(I_\bullet)}$.

For more general subsets $Z\sset \R^m$ and more general filtrations the derivatives/elements of completion are essentially different objects.
  This  case is more involved and the surjectivity of
 completion does not follow from Whitney extension theorem.
In this appendix we prove a sufficient condition for the surjectivity. This ensures the surjectivity for a very broad class of filtrations.

\

Take an open subset $\cU\sseteq\R^m$ and the ring $R= \quots{C^\infty(\cU)}{J}$. Take a filtration $I_\bullet$ of $R$
 and the corresponding completion  $R\to \hR^{(I_\bullet)}$.
 The elements of $\hR^{(I_\bullet)}$ are  Cauchy sequences of functions, $\{f_j\}\in  R $, such that $f_{j+i}-f_j\in I_j$, for all $i,j>0$.
  Equivalently, the elements can be presented as (formal) sums $\sum^\infty_{j=0}g_j$, for $g_j\in I_j$.

\subsection{Surjectivity for $C^\infty(\cU)$ vs surjectivity for $\quots{C^\infty(\cU)}{J}$}
Fix a filtration $I_\bullet$ of $C^\infty(\cU)$. Assume $I_j\supseteq J$ for any $j$, thus we have the induced
 filtration $\{\pi(I_j):=\quots{I_j}{J}\}$ of $\quots{C^\infty(\cU)}{J}$
 and the diagram
\\\parbox{11cm}
{on the right. The maps $\pi$, $\hat{\pi}$ are surjective. Thus the surjectivity of $\phi$ implies that of $\phi_{/J}$.
Vice versa, assume that $\phi_{/J}$ is surjective.
 Fix an element $\hg\in \widehat{C^\infty(\cU)}^{(I_\bullet)}$ and take any element $g\in \pi^{-1}\phi^{-1}_{/J} \hat\pi(\hg)\sseteq S$.
 Then $\phi(g)-\hg\in \phi(J)=0\in \widehat{C^\infty(\cU)}^{(I_\bullet)}$. Thus $\phi$ is surjective.
}\quad\quad\quad
  $\bM C^\infty(\cU)&\stackrel{\phi}{\to}&\widehat{C^\infty(\cU)}^{(I_\bullet)}
  \\\quad\downarrow\pi&&  \downarrow\hat\pi
  \\\quots{C^\infty(\cU)}{J}&\stackrel{\phi_{/J}}{\to}&\widehat{\quots{C^\infty(\cU)}{J}}^{(\pi(I_\bullet))}\eM$

 Therefore it is enough to  verify the surjectivity of the completion map $C^\infty(\cU)\to \widehat{C^\infty(\cU)}^{(I_\bullet)}$.

\subsection{A sufficient condition for the surjectivity $C^\infty(\cU) \stackrel{\phi}{\twoheadrightarrow}\widehat{C^\infty(\cU)}^{(I_\bullet)}$}
\bthe\label{Thm.Completion.Cinfty.surjective} Let $R=C^\infty(\cU)$, for an open subset $\cU\sseteq\R^m$.
 Suppose there exists an open cover $\cU=\cup \cU_\al$ such that,
 when restricted to each $\cU_\al$, the filtration $I_\bullet$ is  equivalent to the filtration $\ca_0+\sum_k \ca_k\cdot \cb_{k,\bullet}$, where
  (all the ideals depend on $\cU_\al$):
\bei
\item The ideals $\ca_0$, $\{\ca_k\}$ do not depend on $\bullet$;  the  collection $\{\ca_k\}$ is finite and $\{\ca_k\}$  are all finitely generated.
\item The zero loci satisfy: $V(\cb_{k,j})=V(\cb_{k,1})$ for any $k,j$.
\item The ideals $\{\cb_{k,j}\}$ satisfy: $\cb_{k,j}\sseteq I(V(\cb_{k,1}))^{d_j}$, for a sequence $d_j\to\infty$.
\eei
Then the following holds:
\bee[1.]
\item The $I_\bullet$-completion map is surjective, $R\twoheadrightarrow \hR^{(I_\bullet)}$.
\item Moreover, if a closed subset $Z\sset \cU$ satisfies $I_\infty\supseteq I(Z)^\infty$ then any element  $\hf\in \hR^{(I_\bullet)}$ admits a preimage
 which is real analytic off $Z$, i.e.  $f\in C^\infty(\cU)\cap C^w(\cU\smin Z)$.
\eee
\ethe
 \bpr
  Given $\sum g_j$, with $g_j\in I_j\sset C^\infty(\cU)$, we should construct a  function  $f\in
C^\infty(\cU)$,  satisfying:
\beq
\forall\ N:\ f-\suml^N_{j=0}g_j\in I_N.
\eeq
First we reduce the proof to the ring $C^\infty(Ball_1(o))$ and a very particular filtration. Then we estimate the growth of derivatives of $g_j$.
 Then we construct $f$ from $\{g_j\}$ using the cutoff functions with controlled growth.
 Finally, in Step 5, we use Whitney approximation theorem to achieve a function real-analytic
  off $Z$.

\bee[\bf Step 1.]
\item  (Simplifying the filtration $I_\bullet$)
We reduce the statement to the particular case of the filtration $I_\bullet$ of $C^\infty(Ball_1(o))$ satisfying:
\beq\label{Eq.particular.filtration}
\{V(I_j)=V(I_1)\}_j \quad
\text{ and }\quad \{I_j\sseteq I(V(I_1))^j\}_j
\eeq
\bee[\bf i.]
\item Take an open cover by small balls, $\cU=\cup Ball_\al $, such that on each ball $I_\bullet$ is equivalent to the corresponding
 $\{\ca_0+\sum_k \ca_k\cdot \cb_{k,j}\}_j$. We can assume that this covering is locally finite (by shrinking the balls if needed).
  Take the corresponding partition of unity,
\beq
\{u_\al\in C^\infty(\cU)\}_\al:\quad\quad
0<u_\al|_{ Ball_ \al} \le1,\quad\quad
 u_\al|_{\cU\smin Ball_\al}=0, \quad\quad
 \sum u_\al=\one_\cU.
\eeq
Suppose we have proved the surjectivity on each ball.
 Thus for any $\sum g_j\in \widehat{C^\infty}(\cU)^{(I_\bullet)}$ and each $Ball_\al$ the element $\sum_j u_\al g_j$ is realized, i.e.
  we have $f_\al\in C^\infty(Ball_\al)$ satisfying:
\beq
\forall\ N:\quad u_\al f_\al-\sum^N_{j=0}u_\al g_j \in I_{N+1}\cdot C^\infty(Ball_\al).
\eeq
Then $f:=\sum u_\al f_\al$ is the needed function. Indeed, $f\in C^\infty(\cU)$, as the sum is locally finite, and
 \beq
 f-\sum_j g_j=\sum_\al u_\al f_\al-\sum_j \one_\cU\cdot g_j=\sum_\al u_\al(f_\al-\sum_j   g_j)\in I_\infty.
 \eeq

\

Thus we restrict to $C^\infty(Ball_1(o))$ and replace $I_\bullet$ by the equivalent filtration as in the
assumptions. Thus
 $I_j=\ca_0+\sum_k \ca_k\cdot \cb_{k,j}$.

\item
  An element of $\hR$ is $c_0+\sum_{j\ge1} (g^{(0)}_j+g^{(>0)}_j)$, where $g^{(0)}_j\in \ca_0$ and
   $g^{(>0)}_j \in \sum_k \ca_k\cdot \cb_{k,j} $.
   We should construct $f\in R$ that satisfies: $f-c_0-\sum^N_{j=1} (g^{(0)}_j-g^{(>0)}_j)\in I_{N+1}$ for any $N$.
   As $\ca^0\sseteq I_j$, for each $j$, one can omit $g^{(0)}_j$.
 This reduces the statement to the filtration $\{\sum_k \ca_k\cdot \cb_{k,j}\}_j$.

\item
Suppose $\{I_j=\sum_k \ca_k\cdot \cb_{k,j}\}_j$. For each $\ca_k$
fix a (finite)  set of generators, $\{a^{(k)}_i\}_i$.
 Then an element $g_j\in I_j$ is written as $\sum_{k,i} a^{(k)}_i\cdot b_{i,k,j}$, with $b_{i,k,j}\in \cb_{k,j}$.
 Therefore  $\sum g_j\in \hR$ is presentable as $\sum_{k,i} a^{(k)}_i \big(\sum_{j=0}^\infty b_{i,k,j}\big)$. (Here the sum over $i,k$ is finite.)
 It is enough to find some $C^\infty$-representatives $\{\tb_{i,k}\}$ of $\{\sum_j b_{i,k,j}\}$, i.e $\tb_{i,k}-\sum_j b_{i,k,j}\in \cap_j \cb_{k,j}$.
  Indeed, for such representatives we get:
 \beq
 \sum_{k,i} a^{(k)}_i \tb_{i,k} -\sum_{k,i} a^{(k)}_i \big(\sum_{j=0}^\infty b_{i,k,j}\big)\in \capl_j\sum \ca_k\cdot\cb_{k,j}=I_\infty.
\eeq Thus it is enough to consider just the filtration $\{\cb_j\}$,
with $\cb_j\sseteq I(V(\cb_1))^{d_j}$, for a sequence
$d_j\to\infty$.
 Furthermore, we pass to an equivalent filtration satisfying  $\cb_j\sseteq I(V(\cb_1))^j$.

\

Therefore it is enough to establish the surjectivity
$R\twoheadrightarrow \hR^{(I_\bullet)}$ for the  filtration of the
particular type
 as in equation \eqref{Eq.particular.filtration}.
 \eee

\item
We have $\{g_j\in I_j\}$ for the specific filtration of the ring $C^\infty(Ball_1(o))$ as in \eqref{Eq.particular.filtration}.
 By slightly shrinking the ball we can assume $g_j\in C^\infty(\overline{Ball_1(o)})$, in particular each derivative of each $g_j$ is bounded.

We claim, for any $j$ and  any $\uk$ with $|\uk|< j$, and any $x\in Ball_1(o)$ the derivatives are bounded:
\beq
  |g^{(\uk)}_j(x)|<C_{g_j}\cdot dist(x,Z)^{j-|\uk|}.
\eeq
(Here $\{C_{g_j}\}$ are some constants that depend on $g_j$ only.)

Indeed, fix some $x\in Ball_1(o)\smin Z$ and some $z\in Z$ for which $dist(x,z)-dist(x,Z)\ll dist(x,Z)$.
 By the assumption $g_j\in \cm_z^j$, thus $g^{(\uk)}_j|_z=0$ for $|\uk|<j$. Therefore the Taylor expansion  with remainder (in $Ball_{dist(x,z)}(z)$)
  gives:
\beq
g_j( x )=\sum_{|\uk|=j}\frac{|\uk|}{\uk!}\Big(\intl^1_0(1-t)^{|\uk|-1}
g^{(\uk)}_j|_{(\uz+t(\ux-\uz))}dt\Big) ( \ux -\uz)^\uk. \eeq (Here
$\ux,\uz$ are the coordinates of $x,z$, $\uk!=k_1!\cdots k_m!$, and
$g^{(\uk)}(\dots)$ is a multi-linear form.)

    Note that $|(\ux -\uz)^\uk|\le \big(dist( x ,Z)+\ep\big)^j$
  and the derivatives   $g^{(\uk)}_j$ are bounded on  $Ball_1(o)$. Thus
  $|g_j( x )|\le C_0\cdot dist( x ,Z)^j$, for a constant $C_0$.

The bounds on the derivatives,  $|g^{(\uk)}( x )|\le\dots$, are obtained
in the same way, by Taylor expanding $g^{(\uk)}$ at $z$.

\item
 We use a particular cutoff function with controlled growth of derivatives:

{\bf Theorem 1.4.2 of \cite[pg. 25]{Hormander}} {\em  For any compact set with its neighborhood, $Z\sset \cU\sset\R^n$,  and
 a positive decreasing sequence  $\{d_j\}$  satisfying $\sum d_j< dist(Z,\di\cU)$, there exists a smaller neighborhood, $Z\sset \cV\ssetneq \cU$,  and
     a    function $\tau\in C^\infty(\R^n)$ satisfying
\bee
\item $\tau|_{\R^n\smin \cU}=0$, $\tau|_{\cV}=1$;
\item for any $\uk$, and $x,y_1,\dots,y_k\in \R^n$ the bound holds: $|\tau^{(\uk)}(x)(y_1,\dots,y_k)|\le \frac{C^{|\uk|}\cdot |y_1|\cdots|y_k|}{d_1\cdots d_k}$.
\eee
}
(Here   the constant $C$ depends only on the dimension $n$.)

\

In our case the subset  $Z\sset  Ball_1(o) $ is closed and we can assume it is compact by shrinking the ball.
  Define the $\ep$-neighborhood,   $\cU_\ep(Z):=\{x|\ dist(x,Z)<\ep\}\sset \R^n$.
 Fix a  decreasing sequence of positive numbers $\{\ep_j\}$,   $\ep_j\to 0$.
  Assume it decreases fast, so that for each $j$ exists a cutoff function satisfying:
 \beq
 \tau_j|_{\cU_{\ep_{j+1}}}=1,\quad\quad
   \tau_j|_{Ball_1(o)\smin \cU_{\ep_j}}=0,\quad\quad
   \text{and $|\tau^{(\uk)}_j|$ is bounded as above, for any $\uk$.}
 \eeq
Define $f(x):=\sum_j  \tau_j(x)\cdot g_j(x)$. We claim that $f\in C^\infty(Ball_1(o))$, when $\{\ep_j\}$ decrease fast.

 The statement $f\in C^\infty(Ball_1(o)\smin Z)$ is obvious, as for any $x\in Ball_1(o)\smin Z$ the summation is finite.
 To check the behaviour on/near $Z$ we bound the derivatives:
\begin{multline}
\Big|\big(\tau_j(x)\cdot g_j(x)\big)^{(\uk)}\Big|=
\Big|\suml_{0\le \ul\le \uk}\bin{|\uk|}{\ul}\tau^{(\ul)}_j(x)\cdot g^{(\uk-\ul)}_j(x)\Big|
<\\<
\suml_{0\le \ul\le \uk}\bin{|\uk|}{\ul}\Big|\tau^{(\ul)}_j(x)\cdot C_{g_j}\cdot dist(x,Z)^{j-|\uk|+|\ul|}\Big| <
\\<\suml_{0\le \ul\le
\uk}\bin{|\uk|}{\ul}C^{|\ul|} \cdot
C_{g_j}\frac{dist(x,Z)^{j-|\uk|+|\ul|}}{d_1\cdots d_{|\ul|}}<
dist(x,Z)\cdot \suml_{0\le \ul\le
\uk}\bin{|\uk|}{\ul}\frac{C^{|\ul|} \cdot C_{g_j} }{d_1\cdots
d_{|\ul|}} \ep_j^{j-|\uk|+|\ul|-1}.
\end{multline}
We assume the sequence $\{\ep_j\}$ decreases fast to ensure, for $j>|\uk|+1$:
\beq
\suml_{0\le \ul\le \uk}\bin{|\uk|}{\ul}\frac{C^{|\ul|} \cdot C_{g_j} }{d_1\cdots d_{|\ul|}}\ep_j^{j-|\uk|+|\ul|-1}<\frac{1}{j!}.
\eeq
Present $f^{(\uk)}(x)=\sum^{|\uk|+1}_{j=0}\big(\tau_j(x)\cdot g_j(x)\big)^{(\uk)}+\sum_{j>|\uk+1|}\dots$.
 Our bounds ensure that the infinite tail converges uniformly on the whole $\R^n$. Thus each $f^{(\uk)}$ is continuous.

\

\item
We claim: $\tau_j \cdot g_j -g_j\in I_\infty$, for any $j$. For this, we construct a function $q\in I_\infty$,   satisfying $Z=q^{-1}(0)$.
 For any $ x_\al\in Ball_1(o)\smin Z$   fix some $q_\al\in I_\infty$ such that $q_\al( x _\al)\neq0$. (This exists as $V(I_\infty)=Z$.)
 By compactness considerations we get a finite subset $\{q_\al\}$ such that the function
   $q( x ):=\sum q^2_\al( x )\in I_\infty$ does not vanish at any point of $Ball_1(o)\smin Z$.

Finally,  $\tau_j \cdot g_j-g_j$ vanishes on $\cU_{\ep_j}$, thus $\frac{\tau_j \cdot g_j-g_j}{q}$ extends to a smooth function on $Ball_1(o)$.
 Therefore $\tau_j \cdot g_j-g_j\in (q)\in I_\infty$.

Hence $f-\sum^N_{j=0} g_j\in I_N$, for any $N$. Thus the completion map sends $f$ to $\sum g_j$.

\item We prove part 2 of the theorem.

Let $I_\bullet$ be a filtration of $C^\infty(\cU)$, as in the assumption. Take an element of the completion, $\sum g_j\in \widehat{C^\infty(\cU)}\ ^{(I_\bullet)}$.
In the previous steps we have constructed a  representative $f\in C^\infty(\cU)$ of $\sum g_j$.
 Take a closed set $Z\sset \cU$ and assume $I_\infty\supseteq I(Z)^\infty$.
  Apply the Whitney extension theorem, see \cite[pg. 65]{Whitney},
   to the restriction of $f$ and all of its derivatives onto $Z$.
 We have the continuous functions $\{f^{(\uk)}|_Z\}_{\uk}$, which  satisfy the compatibility conditions of Whitney.
  (Because they are all restrictions of the derivatives of $f\in C^\infty(\cU)$.)
   Then we get a function $f_{ann}\in C^\infty(\cU)\cap C^w(\cU\smin Z)$,
    whose derivatives (of all orders) on $Z$ coincide with the derivatives of $f$. Which means: $f_{ann}-f\in I(Z)^\infty\sseteq I_\infty$.
    Thus $f_{ann}$ is also a representative of $\sum g_j$.
\epr
\eee

\bex\label{Ex.surjectivity.of.completions}
The class of filtrations of the theorem, locally equivalent to $\big\{\ca_0+\sum_k \ca_k\cdot \cb_{k,j}\big\}$, is rather large.
 In the simplest cases we get various classical results.
 \bee[i.]
\item As the simplest case suppose  $V(I_j)=V(I_1)=V(\cm)$.
 For $I_j=\cm^j\sset \quots{C^\infty(\R^m,o)}{J}$, or more generally when the filtration $I_\bullet$ is equivalent to $ \cm^\bullet$, we get
  the Borel lemma.
\item
Suppose $V(I_j)=V(\cm)$ and $I_j\sseteq \cm^{d_j}$, with $d_j\to
\infty$, but $I_j\not\supseteq\cm^{N_j}$, for any $N_j<\infty$.
 (This happens, e.g. when $I_j$ is generated by flat functions.)
We still get the surjectivity of completion, though not
 implied by Borel's lemma:
 for $\sum g_j\in \widehat{\quots{C^\infty(\R^m,o)}{J}}\ ^{(I_\bullet)}$ we have a representative
   $f\in \quots{C^\infty(\R^m,o)}{J}$, with $f-\sum g_j\in I_\infty$.

\item More generally, suppose $V(I_j)=V(I_1)=:Z$ and $I_j\sseteq (I(Z))^{d_j}$, for $d_j\to\infty$. Again, theorem \ref{Thm.Completion.Cinfty.surjective}
 implies the surjectivity of completion. Note that we do not assume any regularity/subanalyticity conditions on the closed set $Z$.

If $Z\sset \cU$ is a discrete subset then we get a ``multi-Borel" lemma.

\item Take the ring $C^\infty(\R^m_x\times\R^n_y,o)$ with coordinates $x,y$, and the filtration $(y)^\bullet$. The completion map is the
Taylor map in $y$-coordinates, and theorem \ref{Thm.Completion.Cinfty.surjective} ensures its surjectivity:
\[
C^\infty(\R^m\times\R^n,o)\twoheadrightarrow \widehat{C^\infty(\R^m\times\R^n,o)}\ ^{(y)^\bullet}=C^\infty(\R^m,o)[[y]].
\]
(And, moreover, the preimage can be chosen $y$-analytic for $y\neq 0$.)
This recovers the classical Borel  theorem, see   \cite[Theorem 1.2.6, pg. 16]{Hormander} and \cite[Theorem 1.3, pg. 18]{Moerdijk-Reyes}.

\item Many important filtrations are not of the type $I^\bullet$, and not equivalent to this.
 For example, take a hypersurface germ $(X,o)\sset (\R^m,o)$, and assume its singular locus is  $\{x_1=0=x_2\}\sset (\R^m,o)$.
 (This is the usual model case in Singularity Theory, when studying non-isolated singularities.)
 The standard filtration in this case is:
\beq
I_j=(x_1,x_2)^2\cdot(x_1,\dots,x_m)^j\sset C^\infty(\R^m,o).
\eeq
More generally,
 $I_j=x^2_1(x_1,y_1)^{n_{1,j}}+x^2_2(x_2,y_2)^{n_{2,j}}+\cdots +x^2_m(x_m,y_m)^{n_{m,j}}\sset C^\infty(\R^m_x\times \R^n_y,o)$ is
  a typical filtration for complete intersections with non-isolated singularities.

Using the surjectivity of completion, one pulls-back various formal results, over $\R[[x,y]]$, to the $C^\infty$-statements.

   \item
   For the ring $R=\quots{C^\infty\big((\R^m,o)\times(0,1)\big)}{J}$, we can interpret the elements as the families of function germs.
    Then we   get the surjectivity of completion in families,
  $R\twoheadrightarrow \hR$.
 For example, for $I_j=( \ux )^j$ we get a particular version of Borel lemma in families: any power series
$\sum a_{\um}(t) \ux ^\um$, with $ a_{\um}(t)\in C^\infty([0,1])$,
 is the $\ux $-Taylor expansion of some function germ $f_t( x )\in C^\infty\big((\R^m,o)\times(0,1)\big)$.
 \eee
\eex

\beR\label{Ex.surjectivity.of.completions.remark}
 For some rings/filtrations  one can  apply the following surjectivity
argument. Assume $(R,\cm)$ is local and $I_j\sseteq \cm^{d_j}$, with
$d_j\to \infty$.   Then the completion $R\to \hR^{(\cm)}$ factorizes
through
 $R\to \hR^{(I_\bullet)}\to \hR^{(\cm)}$. Thus, if the map $R\to \hR^{(\cm)}$ is surjective and the map  $\hR^{(I_\bullet)}\to\hR^{(\cm)}$
  is injective, the map $R\to \hR^{(I_\bullet)}$ is surjective.
 However, a necessary condition for the injectivity   $\hR^{(I_\bullet)}\twoheadrightarrow\hR^{(\cm)}$ is $I_\infty\supseteq\cm^\infty$.
  And this does not hold for many filtrations of $C^\infty(\cU)$.
\eeR

\end{document}